\documentclass[12pts]{article}
\usepackage{amsmath,amsthm,amssymb,amscd}
\usepackage{latexsym}
\newcommand{\oths}{{\theta}^{\#}}

\newcommand{\talpha}{\tilde{\alpha}}
\newcommand{\tbeta}{\tilde {\beta}}

\newcommand{\oa}{\bar {{\cal A}}}
\newcommand{\da}{d_{\cal A}}

\newcommand{\ot}{\oths d\tau}
\newcommand{\dc}{\delta}
\newcommand{\xwe}{ \xi\wedge \eta }
\newcommand{\xec}{[\xi, \eta ]^{\circ}}
\newcommand{\xe}{[\xi, \eta ]}

\newcommand{\be}{\begin{eqnarray*}}
\newcommand{\ee}{\end{eqnarray*}}

\newcommand{\ga}{\Gamma ({\cal A})}

\newcommand{\cm}{C^{\infty}(M)}
\newcommand{\cdm}{C^{\infty}(M)^D}
\newcommand{\mg}{M\times \f{g}}
\newcommand{\mgg}{M\times \f{g^*}}
\newcommand{\Lam}{\Lambda}

\newcommand{\rt}{\Lambda}

\newcommand{\f}[1]{{\mathfrak #1}}

\newcommand{\plamm}[1]{\frac{\partial }{\partial x_{#1}}}

\newcommand{\lon}{\longrightarrow}

\newcommand{\reals}{{\mathbb R}}

\newcommand{\half}{\frac{1}{2}}

\newcommand{\cinf}{C^{\infty}}

\newcommand{\cala}{{\cal A}}
\newcommand{\gm }{\Gamma }
\newcommand{\alt}{\mbox{Alt}}
\newcommand{\thh}{\rho_* }
\newcommand{\smalcirc}{\mbox{\tiny{$\circ $}}}
\newcommand{\frakg}{{\frak g}}
\newcommand{\backl}{\mathbin{\vrule width1.5ex height.4pt\vrule height1.5ex}}
\newcommand{\per }{\backl }
\newcommand{\deltaa}{\delta^{\bullet} }
\newcommand{\deltab}{\delta}
\newcommand{\chii}{{\frak X}}

\newtheorem{thm}{Theorem}[section]
\newtheorem{lem}[thm]{Lemma}
\newtheorem{cor}[thm]{Corollary}
\newtheorem{pro}[thm]{Proposition}
\newtheorem{ex}[thm]{Example}
\newtheorem{rmk}[thm]{Remark}
\newcommand{\pf}{\noindent{\bf Proof.}\ }
\newtheorem{defi}[thm]{Definition}

\begin{document}

\title{{\bf The local structure of  Lie bialgebroids}}

\author{Zhang-Ju Liu
\thanks{Research partially supported
by NSF of China and the Research Project of ``Nonlinear Science".}\\
 Laboratory of Pure and Applied  Mathematics   \\
Peking  University \\
Beijing, 100871, China\\
        {\sf email: liuzj@pku.edu.cn}\\
        Ping Xu
         \thanks{ Research partially supported by NSF
       grant DMS00-72171. }\\
        Department of Mathematics\\
         Pennsylvania State University \\
         University Park, PA 16802, USA \\
{\sf email: ping@math.psu.edu }}

\date{}
\maketitle \centerline{\bf Dedicated to Professor Qian Min on the
occasion of his 75th birthday}

\begin{abstract}
We study the local structure of Lie bialgebroids at regular
points. In particular, we classify all transitive Lie
bialgebroids. In special cases, they are connected to classical
  dynamical $r$-matrices
and  matched pairs induced by Poisson   group actions.
\end{abstract}

\section{Introduction}

    Lie bialgebroids  were   introduced by
Mackenzie and Xu in \cite{MX} as the infinitesimal versions of
Poisson groupoids.
Lie algebroids are, in a certain sense,   generalized tangent
bundles, while Lie bialgebroids can be considered as
generalizations of both Poisson structures and Lie bialgebras.
It is therefore tempting to classify  the local structure of  Lie
bialgebroids. It is well known that at generic points (called
regular points) a Lie algebroid $\cala$ is   locally  isomorphic
to the following ``standard" one:

\begin{equation}
\label{standard} {\cal A} =  D \oplus (M \times \f{g}),
\end{equation}
where $D \subseteq TM$ is  an integrable  distribution, which is
the image of the anchor, and $M \times \f{g}$ is a bundle of Lie
algebras, namely the isotropy Lie algebras. Here one should think
of $\f{g}$ as a vector space equipped with a family of Lie algebra
structures parameterized by points on the base manifold $M$,
denoted by $\{\f{g}_x |x\in M\}$. The  Lie brackets
  are  constant along  leaves of $D$.
I.e.,   for any $A, B \in \f{g}$
 considered as constant sections\footnote{Throughout the paper,
 we adapt the convention that
elements in $\f{g}$ (or $\f{g}^*$) are  automatically considered
as constant sections on $\cala$ (or $\cala^*$) unless
 otherwise specified.}
  of ${\cal A}$, $[A, B]\in
 \cdm \otimes \f{g} $, where $ \cdm  $ stands for leafwise constant
 functions on $M$.  The anchor of $\cala $ is the projection
$\cala \lon D$.  The bracket between sections of $\cala$ is given
by
\begin{equation}
\label{eq:mixed} [X+A, Y+B]=[X, Y]+(X(B)-Y(A) +[A, B]), \ \ \ \
\forall X, Y\in \gm (D), \ A, B\in C^{\infty}(M, \f{g} ).
\end{equation}
In particular, $\cala$ is a {\em transitive}  Lie algebroid if
$D=TM$.

Consider the extended  vector bundle
\begin{equation}
\label{oa} \oa := TM \oplus (M \times \f{g}) \supseteq
  D \oplus (M \times \f{g})= {\cal A}.
\end{equation}
 The bracket on  $\gm (\cala )$  naturally extends to a bracket on
 sections of  $\oa$. This is  done simply by applying  the
same formula: Equation (\ref{eq:mixed}) to $X +A,\  Y +B\in \chii
(M) \oplus C^{\infty}(M, \f{g} )$. Using the graded Leibniz rule,
one may extend this bracket to sections of $\wedge^* \oa$ as well,
which extends the ordinary Schouten brackets on $\gm (\wedge^*
{\cal A})$. However,  $\oa$ is  in general {\em not}  a Lie
algebroid any more.
 For example, for $A, B \in \f{g},~ X \in \chii (M)$, we have
\begin{equation}
\label{xab} [X, [A,B]] + [B, [X, A]]+[A,[B, X ]] = L_X [A,B],
\end{equation}
which may not be zero unless $[A, B]$ is a constant. This extended
bundle $\oa$  and the  brackets  on $\gm (\wedge^* \oa )$ will be
useful later in the paper in order to describe the Lie algebroid
structure  on $\cala^*$.

For any section $\Lam \in \gm (\wedge ^2 \oa )$,  by $[\cdot ,
\cdot]_{\Lam}$ we
 denote the bracket on $\gm (\oa^* )$ defined by
\begin{equation}
[\phi, \psi ]_\Lam =L_{\Lam^{\#}\phi} \psi-L_{\Lam^{\#}\psi} \phi
-d[ \Lam (\phi , \psi  )] \ \ \ \forall \phi, \psi \in \gm (\oa^*
).
\end{equation}

The purpose of this paper is to classify  all possible Lie
bialgebroids $(\cala , \cala^* )$ for which  $\cala$ is
 of the standard form (\ref{standard}).
We find  that  a Lie bialgebroid structure in this case  can  be
totally characterized by a  pair $(\Lam , \deltab )$, where $\Lam$
is a section of $\wedge^2 \bar{\cala}$,
 and $\delta $ is a bundle
map $M\times   \frakg \lon M \times \wedge^2 \frakg$. Such a pair
   defines the Lie algebroid
 structure on  $\cala^*$. For transitive
Lie bialgebroids, one can furthermore describe such  a Lie
bialgebroid using a  quadruple $( \pi , \theta, \tau , \dc )$
satisfying certain equations (see Theorem \ref{gauge1}), where
$(\pi , \theta , \tau)$ are components of $\Lam$. In this way, one
establishes a one-to-one correspondence between transitive Lie
bialgebroids  and equivalence classes of certain quadruples. As
applications, we recover a   construction, due to Lu \cite{Lu}, of
a Lie algebroid arising from a  matched pair induced by a Poisson
group action. Another special case is connected to dynamical
$r$-matrices coupled with Poisson manifolds, which is  a natural
generalization of the classical  dynamical $r$-matrices of Felder
\cite{Felder}.

{\bf Acknowledgments.} We  are grateful to Vladimir  Drinfel'd for
raising the local structure question.
 We   also    wish to thank
 Kirill Mackenzie for useful discussions and Alan Weinstein
for numerous suggestions and comments.
 In addition to the funding sources mentioned
in the first footnote, the second author would also  like to thank
Peking University for  hospitality while part of this project was
being done.

\section{Regular Lie bialgebroids}

The aim of  this section is  to  give a  general description of a
Lie bialgebroid $(\cala , \cala^* )$ at a regular point. As
mentioned in Introduction,  any    Lie algebroid
 $\cala$ is always   locally   isomorphic to   the ``standard"
one (\ref{standard}), so we will assume that $\cala$ is isomorphic
to this form throughout the paper. Indeed we say that a Lie
bialgebroid $(\cala , \cala^* )$ is {\em regular} if $\cala$ is of
the  standard form (\ref{standard}). It is clear that $D^* $ is
isomorphic to the quotient bundle
 $T^*M/D^{\perp}$.
By $pr$, we denote the projection $T^*M \lon D^*$. By abuse of
notation, we use the same symbol to denote its  extension
$$pr: T^*M \oplus (M \times \f{g}^*) \lon {\cal A^*}
(\cong D^* \oplus (M \times \f{g}^*) ).$$ For any $ \alpha \in
\Omega^1(M)$ (or more generally
 $\alpha \in \gm (T^*M \oplus (M \times \f{g}^*) )$, by $\talpha$ we denote
its corresponding section in  $D^* $ (or ${\cal A^*}$).
 Hence
\begin{equation}
\label{eq:dAf}
 \da f = \widetilde{df},~ \forall f \in \cm.
\end{equation}
   Even though we are mainly interested  in the  local structure,
   our results hold for more general situations as well.  Indeed
 we only need to
  require  that the first leafwise de-Rham cohomology
group $H^1_D(M)$ vanishes.

The main theorem of this section is the following:

\begin{thm}
\label{thm:2.1} Let $(\cala , \cala^* )$ be a Lie bialgebroid,
where $\cala$  is of the form (\ref{standard}): ${\cal A}= D\oplus
(M\times \f{g}) $.  Assume  that  $H^1_D (M)= 0$. Then the
differential  $d_* : \gm (\wedge^* \cala ) \lon \gm (\wedge^{* +1}
\cala ) $
 corresponding to the Lie algebroid structure on ${\cal A^*}$
is of the form:
\begin{equation}
\label{eq:d}
 d_* =  [ \Lam , ~ \cdot ~ ]   +  \deltab,
\end{equation}
where $ \Lam  \in \Gamma(\wedge^2 \oa   )$ such that $ [ \Lam , ~
\Gamma({\cal A})] \subset \Gamma(\wedge^{ 2} {\cal A})$
   and $\dc$ is a bundle  map $M\times \frakg \lon
 M\times \wedge^2 \frakg  $,
which is constant along leaves of $D$. Moreover the pair $(\Lam ,
\deltab )$ is unique up to a gauge transformation: $\Lam\to \Lam
+r_0$, $\deltab \to \deltab -[r_0 , \cdot ]$ for some $r_0 \in
C^{\infty}(M, \f{g} )^{D}$.
\end{thm}
Here, as well as in the sequel, $\deltab$  is considered
as  a linear map: $\gm (\wedge^* \cala )
 \lon \gm (\wedge^{* +1} \cala ) $, which is
a graded derivation  naturally
 extending the bundle map $\deltab : M\times \frakg \lon
 M\times \wedge^2 \frakg  $.

We will divide the proof into several lemmas. By $[\cdot, \cdot]_*
$ and  $\rho_* $, we  denote,
 respectively, the Lie bracket and
the anchor of  the  dual Lie algebroid  ${\cal A^*}$.
Define $K\in \gm (\wedge^2 \oa )$ and $\theta \in \chii (M)
\otimes \f{g} $ by
\begin{eqnarray}
&&K(\alpha +\xi , \beta +\eta )= <\rho_{*}\xi , \ \beta
>-<\rho_{*} \eta ,\
 \alpha>;\\
&&\theta (\alpha , \xi )=- <\rho_* \xi   , \ \alpha> , \ \ \ \
\forall \alpha, \beta \in \Omega^{1}(M), \ \xi, \eta \in \frakg^*
. \label{eq:theta}
\end{eqnarray}
Note that in general $K$ is not a section of $\wedge^2 \cala$.
This is because the image of  $\rho_{*}$ is not necessarily
contained in $D$. By  $\theta^{\#}$, we denote   the  induced
bundle map $T^*M \lon \mg $:
\begin{equation}
<\theta^{\#}(\alpha ), \xi> =\theta (\alpha , \xi ).
\end{equation}

 It is obvious  to see that
\begin{eqnarray}
&&K^{\#}|_{T^* M}=\theta^{\#}: \  T^*M\lon  \f{g}; \label{eq:Ktheta} \\
&&K^{\#}|_{M\times \f{g}^*} =\rho_{*}|_{M\times \f{g}^*} :\
M\times  \f{g}^*\lon
 TM. \label{eq:Krho}
\end{eqnarray}
Let   $\pi\in \Gamma(\wedge^2 TM)$ be
 the Poisson structure  on $M$ induced from the
 Lie bialgebroid \cite{MX}.

\begin{pro}
\label{pi} Under the same hypothesis as in Theorem \ref{thm:2.1},
we have
\begin{enumerate}
\item $\pi \in \Gamma ( \wedge^2 D)$;
\item  $\pi^{\#}=\rho_*|_{D^*} \smalcirc pr$;
\item $\forall \alpha, \beta \in \Omega^1(M)$,
$ [\talpha, \tbeta ]_* = [\alpha, \ \beta { \tilde{]}}_{\pi}$;
\item for any $f\in C^{\infty}(M)$, $d_{*}f =-\pi^{\#}(df ) - \oths (df)$;
\item for any $f\in C^{\infty}(M), \ \xi \in \frakg^*$,
$ [ \da f, \ \xi ]_* =- \da (\thh (\xi)f) + ad^*_{\oths  (df) }
\xi  $;
\item $\forall  \xi \in \f{g^*},\,  X \in \Gamma (D)$,
$[ \thh (\xi), X] \in \Gamma (D)$;
\item $\forall x \in M $,
$\oths (D^{\perp}_x) \subseteq Z(\f{g_x} )$, where $Z(\f{g_x} )$
is the center of the Lie algebra $\f{g_x}$.
\item  $\forall e \in \ga  $,
$ [ K , e ] \in \Gamma (\wedge^2 {\cal A})$;
\item $K \in \Gamma (\wedge^2 {\cal A} ) $ iff
$\oths (D^{\perp}) =0$, or
 $      Im (\rho_*) \subseteq D$.
\end{enumerate}
\end{pro}
 \pf By definition \cite{MX}, $\pi^{\#}=- \rho \circ\rho^*_* $, which
implies that $\pi^{\#} (T^*M )\subseteq   D $.
 Thus, it follows that $ \pi \in \Gamma ( \wedge^2 D)$.
On the other hand, we also have  $\pi^{\#}=\rho_* \smalcirc
\rho^*$, thus $\pi^{\#}= \rho_*|_{D^*} \smalcirc pr$. This proves
(2). Therefore
\begin{equation}
\label{eq:rho*} \rho_* =\pi^{\#} +  \rho_* |_{M \times \f{g}^*},
\end{equation}
from which (4) follows immediately. (3) follows because $\rho^*$
is a  Lie algebroid morphism from $T^* M$ to ${\cal A^*}$
\cite{MX}, where $T^* M$ is equipped with the standard Lie
algebroid structure induced from the Poisson tensor $ \pi$.

As for (5), we have \be
&& [ \da f ,  \ \xi]_*\\
 &=&- L_{d_*f}\xi\\
 &=& \da <\oths (df), \xi> + i_{\oths (df) }\da \xi\\
  &=&- \da (\thh (\xi)f) + ad^*_{\oths (df) } \xi.
\ee Here we have used the fact that $\da \xi\in \gm (\wedge^2
(M\times  \frakg^* )  )$. Thus (5) is proved.

Given any $f \in \cdm$, we have $df \in \Gamma (D^{\perp})$, and
therefore
    $ \da f= 0$ according to Equation (\ref{eq:dAf}).
By (5), we obtain that  $ \da (\thh (\xi)f) =0 $
 and $ ad^*_{\oths (df) } \xi= 0$.

It is obvious that $ ad^*_{\oths  (df) } \xi= 0$ is equivalent to
(7). Now $ \da (\thh (\xi)f) =0 $ is equivalent to $d   (\thh
(\xi)f) \in \gm ( D^{\perp} )$,
 which implies that $X(\thh (\xi) f )=0$ for any $X\in \gm (D)$.
The latter is equivalent to $[X, \  \thh (\xi) ](f )=0$, which
implies that $[X,  \ \thh (\xi) ]\in \gm (D )$ since $f$ is an
arbitrary function in $\cdm$. Thus (6) is proved.

It is easy to see that (6) implies that $[K, X]\in \gm (\wedge^2
\cala ), \ \ \forall X\in \gm (D)$. Next we need to show that $[K,
A]\in \gm (\wedge^2 \cala )$ for any  $A \in \frakg$ considered as
a  constant section of $\cala$. For this purpose,  let us choose a
local coordinate system $(x_1, \cdots , x_k , x_{k+1}, \cdots ,
x_{k+l})$ such that $D $ is spanned by $\{ \plamm{k+j}| \ j=1,
\cdots , l\}$. Then clearly  $dx^{i} \in \gm (D^{\perp}), \ i=1,
\cdots , k$. Let $A_{i} (x) =-\oths  (dx^{i} )$.
 According to  (7),
we know that $A_{i}(x)\in Z(\frakg_x )$, for any $x$ and $i=1,
\cdots , k$.
 Thus,  locally we can write $K=K_{1}+K_{2}$, where $K_1 =\sum_{i=1}^{k}
A_{i}\wedge \plamm{i}$ and $K_2 =\sum_{j=1}^{l} A_{k+j} \wedge
\plamm{k+j}$. It thus follows that $[K, A]=[K_{1}+K_{2}, A]=[K_2 ,
A]\in \gm (\wedge^2 \cala )$. Then (8) follows immediately.

Finally, we see, by definition, that $ K \in \Gamma (\wedge^2
{\cal A} )$ iff
 $\oths (D^{\perp}) =0$ or  $Im (\rho_* |_{M \times \f{g}^*})
\subseteq D$. The latter is equivalent to
 $Im (\rho_*) \subseteq D$  due to  (1) and
Equation (\ref{eq:rho*}).
This completes the proof of the proposition.  \qed\\\\\\

\begin{cor}
\label{cor:K} We have the following
\begin{enumerate}
\item for any  $X, Y\in \Gamma (D)$ and $A\in \f{g}$, we have
\begin{eqnarray}
&&[K, [X,Y]] =  [[K, X], Y]+[X, [K, Y]]; \label{xy}\\
&& [K, [X, A]] =[[K, X], A] + [X, [K, A]]. \label{xa}
\end{eqnarray}
\item For  $A, B\in \f{g}$, if $L_{\rho_*(e)}[A, B] = 0 $ for any
 $e \in \Gamma ({\cal A}^* |_{x})$,   then
\begin{equation}
\label{ab}
 [K, [A, B]]|_{x}=[[K, A], B]|_{x} + [A, [K, B]] |_{x}.
\end{equation}
\end{enumerate}
\end{cor}
\pf Equation (\ref{xy}) is obvious. For Equation  (\ref{xa}),  we
write $K=K_{1}+K_{2}$ as in the proof of Proposition \ref{pi}.
From the fact that $A_{i}(x)\in Z(\frakg_x )$, $i=1, \cdots , k$,
 we have
$$[K_1 , [X, A]] =[[K_1 , X], A] + [X, [K_1 , A]].$$
On the other hand, from the Jacobi identity of the Schouten
bracket on $\gm (\wedge^* \cala )$, it follows that
$$[K_2 , [X, A]] =[[K_2 , X], A] + [X, [K_2 , A]].$$
Equation  (\ref{xa}) thus follows.

Finally, we note that
$$[K_2 , [A, B]] =[[K_2 , A], B] + [A, [K_2 , B]],$$
  while
$$([K_1 , [A, B]] -[[K_1 , A], B] - [A, [K_1 , B]] )|_{x}=
\sum_{i=1}^{k} A_i \wedge L_{\plamm{i}} [A,B]. $$ For any $\xi \in
\frakg^*$, it is simple to see that $\xi \per (\sum_{i=1}^{k} A_i
\wedge L_{ \plamm{i}} [A,B]) =L_{K_1^{\#}(\xi )}[A,B]$. On the
other hand, $\rho_{*}(\xi )=\pi^{\#} (\xi )+K_1^{\#}(\xi
)+K_2^{\#}(\xi )$. Since $\pi^{\#}(\xi ), \ K_2^{\#}(\xi )\in \gm
(D)$, hence $L_{\rho_{*}(\xi )}[A,B]=L_{K_2^{\#}(\xi )}[A,B]$.
The conclusion (2) thus follows immediately. \qed\\\\\\

The following result describes the bracket of mixed terms in $\gm
(\cala^* )$.

\begin{pro}
\label{cross} For any  $ \alpha \in \Omega^1(M)$ and $\xi \in
C^{\infty}(M, \f{g^*})$, we have
\begin{equation}
\label{eq:cross} [ \talpha , \xi ]_* = [ \alpha , \xi
]^{\tilde{}}_{\pi + K} .
\end{equation}
\end{pro}
\pf First, let $f \in \cm$ be an arbitrary function,
 and  $\xi \in \f{g^*}$ an arbitrary element. Then according to
Proposition \ref{pi} (5), together with Equations
(\ref{eq:Ktheta}, \ref{eq:Krho}), we have   $[ {\tilde df}, \xi
]_* =[ df , \xi ]^{\tilde{}}_{K}$. On the other hand, it is easy
to see that  $[ df , \xi ]_{\pi}=0$. Thus, we have $[ {\tilde  df}
, \xi ]_* = [ df , \xi ]^{\tilde{}}_{\pi + K}$. Now Equation
(\ref{eq:cross}) follows immediately from
the anchor properties of  a Lie algebroid.  \qed\\\\\\

The last  step is to analyze the bracket between elements in
$C^{\infty}(M, \f{g^*})$. Note that $  \forall  \xi, \eta \in
C^{\infty}(M, \f{g^*})$ and $ \forall f \in C^{\infty}(M)$, by
Equation (\ref{eq:Krho}),
 one has
$$ [ \xi, f\eta]_* - [\xi, f\eta ]_K =f([ \xi, \eta]_* - [\xi, \eta ]_K).$$
This  implies that $[\xi, \eta]_* - [\xi, \eta ]_K$ depends on
$\xi, \eta$ algebraically. On the other hand, $[\xi, \eta]_* -
[\xi, \eta ]_K$ can be split into two parts:
\begin{equation}
\label{pw} [\xi, \eta]_* - [\xi, \eta ]_K =\Omega ( \xi, \eta) +
[\xi, \eta ]^{\bullet},
\end{equation}
where $\Omega ( \xi, \eta) \in \Gamma (D^* )$ and $ [\xi, \eta
]^{\bullet} \in  C^{\infty}(M, \f{g^*})$. Now $ [\cdot, \cdot
]^{\bullet}$ defines a  fiberwise bracket  on the bundle $M \times
\f{g^*}$, while $\Omega( \xi, \eta)$ corresponds to
  a $ \f{g}\wedge \f{g}$-valued one-form on  the bundle $D$, i.e.,
$\Omega \in \gm (D^* )\otimes (\wedge^2 \f{g} )$. Let
$\Omega^{\#}: D\lon M\times (\wedge^2 \f{g} )$ be its induced
bundle map, and
 $\deltaa: M\times \frakg\lon M\times (\wedge^2 \frakg )$
 the  fiberwise cobracket
corresponding to $ [\cdot, \cdot ]^{\bullet}$.

\begin{lem}
\begin{equation}
\label{g*3} d_*  = [\pi + K , ~\cdot~] - \Omega^{\#}  + \deltaa
\end{equation}
\end{lem}
\pf By Propositions \ref{pi},  \ref{cross} and the fact that $ [
\alpha , \beta ]_{K} = [\xi, \eta]_{\pi}=0$, we have
\begin{equation}
\label{g*2}
  [ \talpha + \xi, \tbeta +\eta]_*=
  [ \alpha + \xi, \beta +\eta]^{\tilde{}}_{\pi + K}+
   \Omega ( \xi, \eta) + [\xi, \eta ]^{\bullet}, \ \ \ \ \forall
 \alpha, \beta\in \Omega^1 (M), \  \xi, \eta \in C^{\infty}(M, \f{g^*}) ,
 \end{equation}
which implies Equation (\ref{g*3}) immediately.

\begin{pro}
\label{tau} There exists a  function $\tau \in  C^{\infty}( M ,
\wedge^2 \f{g})$
 and a fiberwise
bracket $ [\cdot, \cdot ]^{\circ}$ on the bundle $M \times
\f{g^*}$ constant along leaves of $D$ such that
\begin{equation}
\xe_* =\xe_{K + \tau} + \xec, \ \ \ \ \forall \xi , \eta \in
C^{\infty}(M, \f{g^*}).
\end{equation}
\end{pro}
\pf For any $X, Y \in \Gamma(D)$,  from the   compatibility
condition $d_{*}[X, Y]=[d_{*}X, Y]+[X, d_{*}Y]$ and Equations
(\ref{xy}, \ref{g*3}), it follows  that
 $$
 \Omega^{\#} [X, Y] = [\Omega^{\#}X, Y]+[ X, \Omega^{\#} Y].
 $$
 This is equivalent to that, as a   $ \wedge^2 \f{g}$-valued one-form,
$\Omega$ is $d_{D}$-closed.  By assumption, there exists a
function  $\tau \in  C^{\infty}( M ,  \wedge^2  \f{g})$ such that
$\Omega = d_{D} \tau$. Hence $\Omega^{\#} (X) =L_X \tau, ~~\forall
X \in \Gamma (D)$.

For any $X \in \Gamma(D), A \in \f{g}$,  the compatibility
condition $d_{*}[X, A]=[d_{*}X, A]+[X, d_{*}A]$, together with
Equations (\ref{xa}, \ref{g*3}), implies that $-[\Omega^{\#} (X)
,A]+[X, \deltaa A] =0$. On the other hand,
 we have
\be
 && [ \Omega^{\#} (X), A ] - [X,  \deltaa A]\\
 &=&  [ L_X \tau, A]  - L_X \deltaa A\\
 &=&  L_X ( [ \tau , A] - \deltaa A).
\ee
It thus follows that  $[ \tau , A] - \deltaa A$ is constant along
leaves of $D$. Let $\dc A=\deltaa A-[ \tau , A] $. Thus, the  map
\begin{equation}
\label{dc}
  \dc = \deltaa  - [ \tau , ~\cdot~] :
   \f{g} \lon C^{\infty}(M, \wedge^2\f{g})^D
\end{equation}
is  a  cobracket   which induces a fiberwise  bracket
 $[ \cdot, \cdot ]^{\circ}$ on $\mgg$ constant along
leaves of $D$. Moreover it is simple to see, from Equation
(\ref{dc}), that
$$
[\xi , \eta ]^{\bullet} = \xec + ad^*_{\tau^{\#}\xi}\eta -
ad^*_{\tau^{\#}\eta}\xi  , \ \ \ \ \forall \xi , \eta \in
C^{\infty}(M, \f{g^*}).
$$
Thus, from Equation (\ref{pw}), we get \be [\xi , \eta ]_* &=&
[\xi, \eta  ]_K + \xec + ad^*_{\tau^{\#}\xi}\eta -
ad^*_{\tau^{\#}\eta}\xi
+ (d_{D} \tau ) (\xi, \eta)\\
  &=& [\xi, \eta  ]_K + \xec + [\xi, \eta  ]_{\tau}\\
  &=& [\xi, \eta  ]_{K + \tau} + \xec
\ee This concludes the proof of the  proposition.
\qed \\\\\\
{\bf Proof of Theorem \ref{thm:2.1}}: According to Proposition
\ref{pi} (3), Proposition \ref{cross}, and
 Proposition \ref{tau},  we conclude that
for any $\alpha , \beta \in \gm (D^* ), \ \xi , \eta \in
C^{\infty}(M, \f{g}^* )$,
$$[\alpha +\xi , \beta +\eta ]_{*}=
[\alpha +\xi , \beta +\eta ]_{\Lam}+[\xi , \eta ]^{\circ}, $$
where $\Lambda =\pi + K+\tau$. It thus follows that $ d_* =  [
\Lam , ~ \cdot ~ ]   +  \deltab$, where
 $\deltab :M\times \f{g} \lon M\times \wedge^2 \f{g}$ is the
cobracket dual to $[\cdot , \cdot]^{\circ}$.
 Notice that in  Proposition \ref{tau},  $\tau$ and $\dc$
are not unique. They  can differ  by an element $r_0 \in
C^{\infty}( M ,  \wedge^2 \f{g})^{D} $: $\tau_1 = \tau + r_0$ and
$\dc_1 = \dc - [r_0 , \cdot ] $. However,  $ \deltaa =  \dc + [
\tau , ~\cdot~] = \dc_1 + [ \tau_1 , ~\cdot~]$ is always fixed.
This completes the proof of the theorem.
 \qed \\\\\\

\begin{thm}
\label{thm:2.8} Under the same hypothesis as in Theorem
\ref{thm:2.1}, we have
\begin{enumerate}
\item $ \Lam \in \Gamma(\wedge^2 {\cal A })$  iff
$Im (\rho_*) \subseteq D$.
\item  For any $x\in M$, the cobracket $\dc_x: \f{g}_x  \lon
\wedge^2\f{g}_x$ is  a Lie algebra  1-cocycle
 if $L_{\rho_*(e)}[A, B] = 0 $, $\forall A, B \in \f{g},~
 e \in \Gamma ({\cal A}^* |_{x})$.  In particular,
$ \dc_{x} $ is always a 1-cocycle  when $Im (\rho_*|_x ) \subseteq
D_x$.
\end{enumerate}
\end{thm}
\pf   $K$ is the  only term in $\Lam$ which is not
necessarily a  section in
    $\wedge^2 {\cal A}$. It is  thus clear that
$\Lam  \in \Gamma (\wedge^2{\cal A})$
 iff $Im(\rho_* ) \subseteq D $  according to Proposition  \ref{pi} (9).

For any $ A, B \in \f{g}$, by the compatibility condition
$d_{*}[A, B]=[d_{*}A, B]+[A, d_{*}B]$ and Equation (\ref{g*3}),
 we have
$$
\deltaa [A,B] -( [\deltaa A,B]+[A, \deltaa B]) = -[K, [A,B]] +[[K,
A], B]+[A, [K, B]].
$$
It thus follows from Corollary \ref{cor:K} (2)
 that  $\deltaa_x$ is a cocycle if
$  L_{\rho_{*}|_{x} (\xi  )}[A, B] = 0$ , \,$\forall  \xi \in
\f{g^*}$. Finally, note  that  $\dc_x$ being
 a cocycle is equivalent to    $\deltaa_x$ being a cocycle
since their difference is  a coboundary  $[\tau,~\cdot ~]$.
Therefore, the theorem  is proved.
\qed\\\\\\




\begin{defi}
\label{defi:dec} A regular  Lie bialgebroid $(\cala, \cala^*)$,
where $\cala$ is of the standard form: ${\cal A} =  D \oplus (M
\times \f{g}) $, is  called {\em decomposable} if $Im (\rho_*)
\subseteq D$.
\end{defi}

We should note that the role of $\cala$ and $\cala^*$ is not
symmetric here. In other words,    that $(\cala, \cala^*)$ is decomposable
does not necessarily mean that $(\cala^* , \cala )$ is decomposable.
In fact $(\cala^* , \cala )$ may even  not be regular.

The following immediately follows from   Proposition \ref{pi} (9).

\begin{cor}
Given a regular  Lie bialgebroid $(\cala, \cala^*)$, if $\f{g_x},
\forall x \in M$,  are  center
 free (e.g., semisimple), then $(\cala, \cala^*)$
is decomposable.

In particular, if  ${\cal A} = M\times \f{g} $ is  a bundle of Lie
algebras and  $\f{g_x}, \forall x \in M$,  are  center free,
 ${\cal A^*} = M\times \f{g^*} $ is  also a Lie algebra bundle such that
 $(\f{g_x}, \f{g^*_x}), \forall x \in M$, are all  Lie bialgebras.
\end{cor}

The following result is obvious.

\begin{pro}
If $(\cala, \cala^*)$ is  a {\em decomposable}  Lie bialgebroid as
in Definition  \ref{defi:dec} and $L$ is any leaf of $D$, then
$(\cala|_L   , \cala^* |_L   )$ is  a transitive  Lie bialgebroid.
\end{pro}

In other words, one can reduce the study of a  decomposable Lie
bialgebroid  to  the study of transitive ones, which is the main
topic of the next section.

\section{Transitive Lie bialgebroids}

This section is devoted to the study of local structures
of  transitive Lie bialgebroids, which is   a special case of regular
ones.  In this case,   we have $\cala =TM \oplus (M \times \f{g}) $ and the
fiberwise  brackets on $M \times \f{g} \lon M$ are constant. We  also
assume, throughout the section, that  $H^1(M )= \{0\}$.
Theorem \ref{thm:2.1} implies the following:

\begin{thm}
\label{gauge}
Let $\cala =TM \oplus (M \times \f{g}) $ be a transitive Lie
algebroid. Assume that $H^1(M )= \{0\}$. Then there is a
one-to-one correspondence between  Lie bialgebroids
 $(\cala , \cala^* )$ and equivalence  classes $[(\Lam , \dc )]$,
 where $ \Lam \in \Gamma(\wedge^2 {\cal A})$,
$ \deltab: \f{g} \lon \f{g}\wedge \f{g} $ is a Lie algebra
1-cocycle satisfying the property that
\begin{equation}
\label{eq:main} [\deltab \Lam +\half [\Lam , \Lam  ], \
e]+\deltab^2 e=0, \ \ \ \forall e\in \gm (\cala ),
\end{equation}
and the equivalence  stands for  the 
gauge equivalence: $\Lam \to \Lam +r_0$, $\dc\to \dc-[r_0 , \cdot
]$ for some $r_0 \in \wedge^2 \f{g}$.
\end{thm}
\pf  Assume that $(\cala , \cala^* )$ is a transitive Lie
bialgebroid.
 According to  Theorem \ref{thm:2.1},  we know that
\begin{equation}
\label{eq:d1}
 d_* =  [ \Lam , ~ \cdot ~ ]   +  \deltab
\end{equation}
for  $ \Lam  \in \Gamma(\wedge^2 \cala )$ and $ \deltab: \f{g}
\lon \f{g}\wedge \f{g} $. Moreover, by Theorem \ref{thm:2.8},  $
\deltab$  is a Lie algebra 1-cocycle. It is simple to check that
\begin{equation}
\label{eq:d2}
 d_{*}^2 =[\deltab \Lam +\half [\Lam , \Lam  ],  ~ \cdot  ]+\deltab^2.
\end{equation}
Thus Equation (\ref{eq:main}) follows.

Conversely, given a pair $(\Lam , \deltab )$ satisfying Equation
(\ref{eq:main}), let $d_* : \gm (\wedge^* \cala )\lon \gm
(\wedge^{*+1} \cala ) $ be as in  Equation (\ref{eq:d1}). Then
$d_* $ defines a Lie algebroid on $\cala^*$ iff $d_{*}^2 =0$,
which is equivalent to: $d_{*}^2 f=0 $ and $d_*^2 e=0$ for any
$f\in C^{\infty}(M)$ and $e\in \gm (\cala )$. Now we easily see
that $\forall f\in C^{\infty}(M), \ e\in \gm (\cala )$,
$$d_*^2  (fe)=fd_*^2  e+(d_*^2 f)\wedge e. $$
It thus follows that $d_*^2 f =0$  whenever the rank of $\cala$ is
greater than $2$. If the rank of $\cala$ is less than or equal to
$2$, by Equation (\ref{eq:d2}), we have $d_*^2 f=[\deltab \Lam
+\half [\Lam , \Lam  ], \ ~ f ~ ]=0$ since $\deltab \Lam +\half
[\Lam , \Lam  ]$ vanishes automatically. Finally, it is clear that
the compatibility condition is satisfied
 automatically.
 \qed\\\\\\


Now  $\Lam$ can be  split   into three parts: $ \Lam =
 \pi + K + \tau $,
where $\pi \in \gm (\wedge^2 TM)$, $K\in \gm (\frakg \wedge TM)$
 and $\tau \in C^{\infty}(M, \wedge^2\f{g})$. Our next task is
 to spell out the
meaning of Equation  (\ref{eq:main}) in terms of these data. Let
us explain some notations that will be  needed below. Recall the
element $\theta$ of $ \chii (M) \otimes \f{g}$
 defined by  Equation (\ref{eq:theta}), and the corresponding
bundle map $\oths: T^*M \lon \mg $.
Note that $\deltab$ extends naturally to a map, denoted by the
same symbol, from $\chii (M)\otimes \f{g}$ to $\chii (M)\otimes
(\wedge^2 \f{g})$.
 If $\theta = \sum X_i \otimes A_i  \in  \chii (M)  \otimes \f{g}$, write
\begin{eqnarray}
\label{dtt} &&\deltab \theta = \sum X_i \otimes   \deltab A_i ,
~~~~~~~~~~~~~~~~~~~~~
[\tau, \theta ]=\sum X_i \otimes [\tau,  A_i ]\\
&&[\theta , \theta] = \sum [X_i, X_j] \otimes (A_i \wedge A_j)
,~~~~ \theta\wedge \theta =\sum X_i \wedge  X_j \otimes [A_i , A_j
].
\end{eqnarray}

Note that $ \ot \in \cinf (M, \wedge^2\f{g}  \otimes \f{g} )$. By
$ \alt ({\ot}  ) \in \cinf (M, \wedge^3\f{g})$,
 we denote its total anti-symmetrization:
\begin{equation}
< \alt ({\ot}  ), \ \xwe \wedge \zeta > = < (\ot )(\xi,\eta),
\zeta > + c. p., \ \ \ \forall \xi, \eta , \zeta \in \frakg^*.
\end{equation}

>From the above discussion, a transitive   Lie bialgebroid
is then  determined by a quadruple $( \pi , \theta, \tau , \dc )$.
We say two quadruples $( \pi_i   , \theta_i , \tau_i  , \dc_i )$,
$i=1, 2$, are  {\em equivalent} if $\pi_1=\pi_2$, $\theta_1=
\theta_2$, and  $\tau_1=\tau_2 +r_{0}$ and $\dc_1=\dc_2 -[r_{0},
\cdot  ]$ for some $r_0\in \wedge^2\frakg$.

We are now ready to state the main theorem of this section.

\begin{thm}
\label{gauge1}
Under the same hypothesis as in Theorem \ref{gauge},
there is a one-one
correspondence between  transitive Lie bialgebroids
$(\cala , \cala^* )$ and
equivalence classes of  quadruples  $( \pi , \theta, \tau , \dc )$
satisfying  the following properties:
\begin{enumerate}
\item $ \dc: \f{g} \lon \f{g}\wedge \f{g} $ is a Lie algebra 1-cocycle;
\item $\pi \in \gm (\wedge^2 TM)$ is a Poisson tensor;
\item $\oths : (T^*M, \pi) \lon \f{g}$ is Lie algebroid morphism;
\item $\deltab \theta + \half [\theta, \theta] =[\tau, \theta]- \pi^{\#}(d\tau) $;
\item $\dc \tau +\half [\tau , \tau] + \alt({\ot})\in \wedge^3 \f{g}$
  is constant on $M$;
\item $\dc^2  + [\dc \tau +\half [\tau , \tau] + \alt({\ot}), \ \cdot ] = 0$, as a map from  $\f{g}$ to  $  \wedge^3\f{g}$.
\end{enumerate}
\end{thm}
\pf Each fiber of $\cala$ is a vector space direct sum
$T_{m}M\oplus \frakg$, and  therefore admits a natural bigrading:
elements in $T_{m}M$ have the degree $(1, 0)$ while elements in
the second component $\frakg$ have the degree $(0, 1)$. This also
induces a bigrading on the fibers of exterior powers of $\cala$
and consequently on their sections.  It is simple to see that
$[\pi , \pi]$ is of degree $(3, 0)$, $[\pi, K]$ is of  degree $(2,
1)$, $[\pi , \tau ]$ is of  degree $(1, 2)$, and $[\tau , \tau ]$
is of  degree $(0, 3)$. On the other hand, $[K, \tau ]$ consists
of elements of degree $(1 ,2)$ and of  $(0, 3)$, and $[K, K]$
consists of elements of degree $(1 ,2)$ and of $(2, 1)$. For any
$S\in \gm (\wedge^3 \cala )$, let $S=\sum_{0\leq i, j\leq 3}S^{(i,
j)}$ be
 its decomposition with respect to this bigrading.
The following lemma can be easily verified by a  direct
computation.

\begin{lem}
\label{lem:3.3} With the above notations,
\begin{enumerate}
\item as a $\wedge^2\frakg$-valued vector field on $M$, we have
$[K, K]^{(1, 2)}=[\theta , \theta ]$;
\item As a $\frakg$-valued bivector field on $M$,
$[K, K]^{(2, 1)} =2 \theta \wedge \theta $;
\item  $\deltab K=\deltab \theta $;
 $[K, \tau ]^{(1, 2)}=-[\tau , \theta ]$;
\item $[K, \tau ]^{(0, 3)}= \alt ({\ot})$;
\item $[\pi , K]=d_{\pi }\theta$.
\end{enumerate}
\end{lem}

Write
$$T=\deltab \Lam +\half [\Lam , \Lam  ] \in \gm (\wedge^3 \cala ) .$$
A direct computation, using Lemma \ref{lem:3.3}, yields that \be
T^{(3, 0)}&=&\half [\pi, \pi ];\\
T^{(2, 1)}&=& [\pi , K]+\half [K, K]^{(2, 1)} \\
&=&\theta \wedge \theta +d_{\pi } \theta; \\
T^{(1, 2)} &=&
 \dc K+ [K, \tau ]^{(1, 2)}+[\pi , \tau ]
+\half [K, K]^{(1, 2)} \\
&=&\deltab \theta -[\tau , \theta ] +\half [\theta , \theta ]+\pi^{\#}(d\tau );\\
T^{(0, 3)}&=&
 \dc \tau +\half [\tau , \tau ] + [K, \tau ]^{(0, 3)} \\
&=&\dc \tau +\half [\tau , \tau ] +\alt ({\ot} ). \ee

For any $X\in {\frak X}(M)$, according to  Equation (\ref{eq:d2}),
we have $d_{*}^2 X=L_{X}T$. Hence  $d_{*}^2 X =0$, $\forall X\in
\chii (M)$,
 is equivalent to that   $T^{(i,  3-i)}
=0$, $i=1, \cdots , 3$,  and $T^{(0, 3)}$ is a constant function.
The latter is   equivalent to       Conditions (2)-(5).

For any $A\in C^{\infty}(M,  \frakg ) $, it is easy to see that
$(d_{*}^2 A)^{(3, 0)}$, $(d_{*}^2 A)^{(2, 1)}$, and $(d_{*}^2
A)^{(1, 2)}$ all vanish by  Conditions (2)-(5). And the only
nontrivial term  remaining is
$$(d_{*}^2 A)^{(0, 3)} =\deltab^2 A+[\deltab \tau +
\half [\tau, \tau ]+\alt ({\ot}), \  A ].$$ Therefore, we see that
$d_{*}^2 e=0, \ \forall e\in \gm (\cala )$ iff Conditions (2)-(6)
hold. The conclusion thus follows by
Theorem \ref{gauge}. \qed\\\\\\

\section{Applications}

As applications, in this section, we will consider two special
cases of transitive Lie bialgebroids.  They are connected,
respectively,
 to the Lie algebroids
arising from  Poisson  group actions \cite{Lu}, and dynamical
$r$-matrices \cite{Felder}.

Let $(\cala , \cala^* )$ be a transitive Lie bialgebroid
corresponding to a  quadruple  $( \pi , \theta, \tau , \dc )$ as
in Theorem \ref{gauge1}. From Proposition  \ref{tau},   it is
simple to
  see that the subbundle  $M \times \f{g}^*$
of ${\cal A^*} =  T^* M \oplus (M \times \f{g}^*)$ is a Lie
subalgebroid iff  $d\tau = 0$, or $\tau\in \wedge^2 \f{g}$ is a
constant.

>From now on, we will assume that  $M \times \f{g}^*$ is a Lie subalgebroid.
Thus one may simply  assume that $\tau=0$ via a  gauge
transformation. By Theorem \ref{gauge1} (6), we have $ \dc^2 = 0
$. This implies that $\frakg^*$ is a Lie algebra  such that
$(\f{g} , \f{g^*})$ is indeed  a  Lie bialgebra.  Theorem
\ref{gauge1} (4) is equivalent to that the linear  map $-(\oths
)^* : \frakg^* \lon  {\frak X}(M)$ obtained by taking the opposite
of the   dual  of $\oths$ is a Lie algebra morphism. Hence, it
defines a  right   action of  $\frakg^*$ on $M$.
  Theorem \ref{gauge1} (3)
is equivalent  to saying that this is a Poisson action
(Proposition 5.4
 in  \cite{Xu}).
In conclusion, when $\tau =0$, the conditions  in Theorem
\ref{gauge1} reduce to the statement
 that $( \f{g},  \f{g^*})$  is  a Lie bialgebra  and
 $M$  is  a Poisson  $\frakg^*$-manifold.
In this case, the Lie algebroid structure on ${\cal A^*} =  T^* M
\oplus (M \times \f{g}^*)$ can be described more explicitly.
 First, we note
that the subbundle  $T^*M$ is also a Lie subalgebroid of ${\cal
A^*}$. In other words, both summands of ${\cal A^*}$ are Lie
subalgebroids. Therefore, in order to describe the bracket on $\gm
({\cal A^*})$, it suffices to consider the bracket between the
mixed terms. For this purpose, consider the maps: \be &&\phi:  T^*
M  \lon     CDO (M \times \f{g}^*),~~~~~~
\phi ({\alpha} )(\xi ) =  L_{\pi^{\#}\alpha}\xi + ad^*_{\oths \alpha}\xi  \ \ \ \mbox{and}  \\
&& \psi : M \times \f{g}^* \lon CDO(T^*M),~~~~~~ \psi({\xi} )(
\alpha ) =-< \xi, d(\oths \alpha )>- i_{(\oths )^* \xi}d\alpha  ,
\ee $ \forall  \alpha \in \Omega^1(M)$ and
  $\xi \in C^{\infty}(M, \f{g^*})$.
Here $CDO$ stands for {\em covariant differential operators}
\cite{Mackenzie}, and $<\cdot , \cdot >$  means the pairing
between $\f{g}^* $ and $\f{g}$.
  Now a simple computation, using Proposition \ref{cross},  yields that
\be
&&  [\alpha, \xi]_* \\
 &=&   [\alpha, \xi]_{\pi +K}\\
&=&L_{\pi^{\#} \alpha} \xi +L_{K^{\#}\alpha } \xi -
L_{K^{\#}\xi}\alpha
-d<K^{\#} \alpha , \xi >\\
&=& L_{\pi^{\#} \alpha} \xi +ad^*_{\oths (\alpha )}\xi-< \oths
(\alpha ), d\xi>
+i_{(\oths )^*\xi}d\alpha +d<\xi ,  \oths (\alpha )>\\
&=&  \phi ({\alpha} )(\xi ) -\psi({\xi} )( \alpha ). \ee

As a consequence, we conclude that (1). both $\phi$ and $\psi$
 are Lie algebroid representations,  (2). $(T^* M, M \times \f{g}^* )$
is  a matched pair,  and (3).
 ${\cal A^*} $ is isomorphic to the corresponding
Lie algebroid $ T^* M \bowtie  (M \times \f{g}^*)$ \cite{M:D43}.
This  Lie algebroid was studied in detail by Lu in \cite{Lu}.

Now we can summarize  the discussion above in the following:

\begin{thm}
\label{thm:match} A quadruple  $( \pi , \theta, 0 , \dc )$ as in
Theorem \ref{gauge1} satisfies the conditions in Theorem
\ref{gauge1} iff $\pi $ is a Poisson tensor, $\dc$ defines a Lie
bialgebra and $\theta$ induces a  right $\frakg^*$-Poisson action.
In this case, the  dual Lie algebroid $\cala^*$ is
 isomorphic to  the matched pair $ T^* M \bowtie  (M \times \f{g}^*)$ as mentioned above.

In other words, under the bijection of Theorem \ref{gauge1}, the
classes of  quadruples  $( \pi , \theta, 0 , \dc )$ correspond to
those transitive Lie bialgebroids $(\cala, \cala^* )$ where $\cala^*
$ is the double of
 the matched pairs $(T^* M , \ M \times \f{g}^* )$
of a Poisson group action.
\end{thm}

\begin{rmk}
{\em More generally, instead of $\tau =0$, one may consider
 the situation where Conditions  (4)-(6)  of  Theorem \ref{gauge1}
are replaced by
\begin{eqnarray}
&&\deltab \theta + \half [\theta, \theta] =0; \label{eq:1}\\
&&[\tau, \theta]- \pi^{\#}(d\tau)=0; \label{eq:2} \\
&&\dc^2 =0;  \ \ \ \mbox{ and }\label{eq:3}\\
&&\dc \tau +\half [\tau , \tau] + \alt({\ot}) =0. \label{eq:4}
\end{eqnarray}
Equations (\ref{eq:1}, \ref{eq:3}) imply that $\dc$ defines a Lie
bialgebra and $\theta$ induces a  right $\frakg^*$-Poisson action,
and therefore one can form a Lie algebroid  $T^* M \bowtie  (M
\times \f{g}^*)$ as in  Theorem \ref{thm:match}.  It is simple  to
see that Equations (\ref{eq:2}, \ref{eq:4}) mean that $\tau\in
C^{\infty}(M, \wedge^2 \f{g})$, considered as a section of
$\wedge^2 \cala $, is a Hamiltonian operator,
 i.e., its graph $\gm_{\tau}$ is a Dirac structure of  $\cala \oplus \cala^*$
\cite{LiuWX:1997}. In this case  $\cala^*$ is isomorphic to the
corresponding Lie algebroid $\gm_{\tau}$. This  is a
generalization of  the situation studied in \cite{LX2}. Note that
Equation (\ref{eq:2}) is equivalent to that
\begin{equation}
\label{eq:5} X_{f}\tau =ad_{\oths (df )} \tau, \ \ \ \forall f\in
C^{\infty}(M),
\end{equation}
which implies that $\tau $   is completely  determined, on a
symplectic leaf,
 by its value at a  particular point of the leaf.  It would
be interesting to investigate   what kind of   equation $\tau$
satisfies when being
  viewed as a function on the leaf space (compare with Theorem 3.14 in
 \cite{EV}).}
\end{rmk}

Another special case of Theorem \ref{gauge1} is when the
quadruple $( \pi , \theta, \tau , \dc )$ is equivalent to $(\pi ,
\theta, \tau , 0 )$. In this case, the corresponding Lie
bialgebroid must be a coboundary Lie bialgebroid with $\Lam =\pi
+K+\tau $.


\begin{thm}
\label{thm:dy}
There is a one-one correspondence between 
coboundary  Lie bialgebroids $(\cala , \cala^* )$, where $\cala
=TM \oplus (M \times \f{g}) $,
 and triples  $( \pi , \theta, \tau ) $
satisfying  the  properties:
\begin{enumerate}
\item $\pi \in \gm (\wedge^2 TM)$ is a Poisson tensor;
\item $\oths : (T^*M, \pi) \lon \f{g}$ is Lie algebroid morphism;
\item $\half [\theta, \theta] =[\tau,  \theta ] -\pi^{\#}(d\tau)$;
\item $ \alt({\ot})+\half [\tau , \tau] \in (\wedge^3 \frakg )^{\frakg} $ (i.e., it   is   constant on $M$ as well as  $ad$-invariant).
\end{enumerate}
\end{thm}

\begin{defi}
For a Poisson manifold $(M, \pi)$ and  a Lie algebra $\f{g}$,
assume that  there exists a  tensor $\theta  \in \chii
(M)\otimes\f{g} $ such that $\oths : (T^*M, \pi) \lon \f{g}$ is a
Lie algebroid morphism. A function  $ \tau  \in  C^{\infty}( M ,
\wedge^2  \f{g})$ is called a  dynamical $r$-matrix coupled with
the Poisson manifold $(M, \pi)$
 via $\theta$ if both conditions (3) and (4) in Theorem \ref{thm:dy}
 are satisfied.
 Here $\theta$ is called a coupling tensor, and the equation
\begin{equation}
\label{eq:dy} \alt({\ot})+\half [\tau , \tau]=\Omega \in (\wedge^3
\f{g})^{\f{g} },
\end{equation}
 is called the  generalized {\bf DYBE} coupled with $( \pi , \theta)$.
\end{defi}


With this definition, Theorem \ref{thm:dy}  can be rephrased as
follows.

\begin{thm}
A function $\tau  \in  C^{\infty}( M ,  \wedge^2  \f{g})$ is a
dynamical r-matrix coupled with a   Poisson tensor
 $ \pi \in \gm (\wedge^2 TM )$
 via $\theta$ iff $\rt: = \pi + K + \tau$
defines a coboundary Lie bialgebroid structure for the Lie
algebroid $TM \oplus (M \times \f{g})$.
\end{thm}

 Now we consider an example, which is in fact the main motivation
for  the  above definition.

\begin{ex}
{\em Assume that  $M=\f{h^*}$ where $\f{h} \subseteq \f{g}$ is  a
Lie subalgebra. Then $M$ is  a   Poisson manifold with the
Lie-Poisson structure $\pi$. In this case, $T^{*}M \cong
\f{h}^{*}\times \f{h}$. Let $\oths :T^{*}M\lon \f{g}$ be the
projection: $(\xi , v)\lon v, \forall (\xi  , v)\in \f{h^*} \times
\f{h}$. Clearly $\oths$ is a Lie algebroid morphism. Let us fix a
basis of $\f{h}$, say $\{e_{1}, e_{2}, \cdots ,e_{k}\}$, and let
$(\lambda_{1}, \cdots ,\lambda_{k})$ be its corresponding
coordinate system on $\f{h}^{*}$. Then  we have
 $\theta =\sum_{i} \frac{\partial}{\partial \lambda_{i}}  \otimes e_{i}$,
where $\frac{\partial}{\partial \lambda_{i}}$, $1\leq i\leq k$,
are considered as constant vector fields on $\f{h}^*$. Clearly
$[\theta , \theta ]=0$. Thus, for any $\tau \in
C^{\infty}(\f{h}^{*}, \wedge^{2}\f{g})$, Condition (3) in Theorem
\ref{thm:dy} takes the form $[\tau,  \theta ]   = \pi^{\#}(d\tau)
$, which, according to Equation (\ref{eq:5}), is equivalent to
that
\begin{equation}
\label{couple1}
 ad_{\xi} \tau  =  X_{l_{\xi}}(\tau ) , ~~~~~~~
 \forall \xi \in \f{h}.
\end{equation}
Here $X_{l_{\xi}} $ denotes the Hamiltonian vector field of the
linear function $l_{\xi}$ on $ \f{h^*}$.
 Thus, If $H$ denotes a connected Lie
group with  Lie algebra $ \f{h}$,  Equation (\ref{couple1}) is
equivalent to
 that the map $\tau: \f{h}^{*}\lon \wedge^{2}\f{g}$ is $H$-equivariant.
 Condition (4)  in Theorem \ref{thm:dy}  becomes
 $$ Alt(d\tau) + \half [\tau, \  \tau ] ~\in ~ (\wedge^3 \f{g})^{\f{g}}.$$
In other words, $\tau$ is a  classical dynamical $r$-matrix in the
sense of Felder \cite{Felder, EV}. We thus recover the main result
in \cite{BK-S}.}
\end{ex}

We end the paper with the following

\begin{ex}
{\em Let $(M ,\pi)$ be any Poisson manifold  and $\f{g}= sl(2,
\reals )$ with the standard  generators $\{H, E_+ , E_-\}$:
$$[H, E_{+}]=E_{+}, \ \ [H, E_{-}]=-E_{-}, \ \ [E_{+}, E_{-}]=2H.$$

Take $\theta = X_{f} \otimes  H $, where $f$ is a smooth function
on $M$ and $X_{f} $ its Hamiltonian vector field. Clearly
 $\oths : (T^*M, \pi) \lon \f{g}$ is Lie algebroid morphism
 because $d_{\pi}\theta +\theta \wedge \theta =0$.

Let $$ \tau =  e^{f} H \wedge E_{+}    +     e^{-f} H \wedge E_- +
 E_+ \wedge E_-  ,
$$
which can be considered as a twist of standard r-matrix on $sl(2,
\reals )$. Obviously $[\theta, \theta] =  0 $  and
$$ [\tau ,  \theta ]  =\pi^{\#}(dr) = X_{f} \otimes
(e^{f} H \wedge E_{+}    -    e^{-f} H \wedge E_{-})
$$
so that the condition in Theorem \ref{thm:dy} (3)  holds. For
Theorem \ref{thm:dy} (4), note that $ \alt({\ot}) = 0$ and  $
[\tau , \tau ] = 3 H \wedge E_{+} \wedge E_{-}\in (\wedge^3 \f{g}
)^{\f{g}} $. Hence $\tau$ is indeed a dynamical $r$-matrix coupled
with $\pi$ via $\theta$.

In particular, let  $M= \reals^2$ be equipped with the standard
symplectic structure $\omega =dx\wedge dy$ and $f =ax - by $. Then
we have \be
\pi&=&\frac{\partial}{\partial x}\wedge \frac{\partial}{\partial y}\\
\theta& =& (b    {\partial \over \partial x}+a {\partial \over
\partial y})
\otimes H,  \ \ \ \mbox{ and}\\
\tau (x, y) &= &e^{ax-by}H\wedge E_{+} +e^{-(ax-by)} H\wedge E_-
+E_{+}\wedge E_{-}. \ee
 }
\end{ex}


\begin{thebibliography}{99}
\bibitem{BK-S}
Bangoura, M. and Kosmann-Schwarzbach, Y., Equation de Yang-Baxter
dynamique classique et algebroides de Lie, {\it C. R. Acad. Sci.
Paris Serie 1}, {\bf 327} 541- 546, (1998).

\bibitem{EV}
  Etingof,  P. and Varchenko A.,  Geometry and
classification of solutions of the classical dynamical Yang-Baxter
equation, {\it Comm. Math. Phys.}, {\bf 192} 77-120, (1998).

\bibitem{Felder}  Felder,  G.,
Conformal field theory and integrable systems associated to
elliptic curves, {\it Proc. ICM  Z\"urich}, Birkh\"auser, Basel,
(1994), 1247-1255.

\bibitem{LiuWX:1997}
 Liu, Z.-J.,  Weinstein, A., and Xu, P.,
 Manin triples for {L}ie bialgebroids,
 {\em J. Differential Geom.} {\bf 45} (1997),  547--574.



\bibitem{LX2}
Liu, Z.-J.,  Xu, P.,
 Dirac structures and dynamical $r$-matrices,
 {\it Ann. Inst. Fourier} {\bf 51} (2001), 831-859.


\bibitem{Lu} Lu, J.-H.,
 Poisson homogeneous spaces and Lie algebroids associated to
  Poisson actions,
 {\em Duke Math. J.}  {\bf 86} (1997), 261--304.




\bibitem{Mackenzie}
Mackenzie, K.,
 Lie groupoids and {L}ie algebroids in differential geometry,
 London Mathematical Society Lecture Note Series, {\bf 124} Cambridge
  University Press, 1987.



\bibitem{M:D43}
 Mackenzie, K.,
 Notions of double for {L}ie algebroids, preprint.

\bibitem{MX}
Mackenzie, K. and Xu, P., Lie bialgebroids and Poisson groupoids,
{\em Duke Math. J.} {\bf 18} (1994), 415-452.

\bibitem{Xu}
Xu, P., On  Poisson groupoids, {\em Internat. J. Math.} {\bf 6}
(1995), 101-124.
\end{thebibliography}
\end{document}